\documentclass[a4paper,12pt,leqno]{article}
\usepackage{latexsym}
\usepackage[all]{xy}
\usepackage{amssymb} 
\usepackage{amsmath} 
\usepackage{theorem}
\usepackage{colortbl} 
\usepackage{graphicx}

\def\Z{{\mathbb{Z}}}
\def\R{{\mathbb{R}}}
\def\A{{\mathcal{A}}}
\def\B{{\mathcal{B}}}
\def\calS{{\mathcal{S}}}
\def\calC{{\mathcal{C}}}
\def\Shi{{\mathcal{S}}}
\newcommand{\bfz}{{\bf z}}

\DeclareMathOperator{\codim}{codim}

\DeclareMathOperator{\Der}{Der}

\numberwithin{equation}{section}
\newcommand{\owari}{\hfill$\square$}

\theoremstyle{break}
\newtheorem{theorem}{Theorem}[section]
\newtheorem{prop}[theorem]{Proposition}
\newtheorem{cor}[theorem]{Corollary}
\newtheorem{corollary}[theorem]{Corollary}

\newtheorem{define}[theorem]{Definition}
\newtheorem{definition}[theorem]{Definition}

\newtheorem{remark}[theorem]{Remark}

\newtheorem{problem}[theorem]{Problem}

\title{Simple-root bases for Shi arrangements}
\author{Takuro Abe 
and Hiroaki Terao
}
\date{\today} 

\begin{document}

\maketitle

\begin{abstract}
In his affirmative answer to the Edelman-Reiner conjecture,
Yoshinaga proved that the logarithmic derivation modules of the cones
of the extended Shi arrangements are free modules.
However, all we know about the bases is their existence
and degrees.
In this article, 
we introduce 
two distinguished bases for the modules.
More specifically, we will define and study
the simple-root basis plus (SRB$+$)
and the simple-root basis minus (SRB$-$)
when a primitive derivation
is fixed.   
They have remarkable 
properties relevant to the simple roots
and those properties 
characterize the bases.
\end{abstract}

\section{Introduction}

Let $V$ be an $\ell$-dimensional Euclidean space.
Let $\Phi$ be an {\bf
irreducible 
(crystallographic)
 root system}
in the dual space $V^{*}$.
Fix
a set
$\Phi^{+} $ 
 of {\bf positive roots}.
For any $\alpha\in\Phi^{+} $ and
 $j\in\Z$,
the affine hyperplane
$$
H_{\alpha, j} 
:=
\{x\in V \mid \alpha(x)=j\}
$$
is a parallel translation of $H_{\alpha} :=H_{\alpha, 0}$. 
The arrangement
$\A(\Phi) := \{H_{\alpha} \mid \alpha\in \Phi^{+} \}
$
is called
the
{\bf
crystallographic arrangement} 
of the type
$\Phi$.

\begin{define}
Let $k\in\Z_{>0} $. 
An \textbf{extended Shi arrangement} $\mathit{Shi}^{k}$ 
of the type $\Phi$ is an affine arrangement defined by 
$$
\mathit{Shi}^{k}
:=
\mathit{Shi}^{k}(\Phi^{+})
:=
\{H_{\alpha, j} ~\mid~ \alpha\in \Phi^{+},\,\,\,j\in\Z, \,-k+1 \leq j\leq k\}.
$$
\label{shi}
\end{define}
The extended Shi arrangements 
for $k=1$ were introduced by 
J.-Y. Shi \cite{Shi1, Shi2}
in his study of the Kazhdan-Lusztig 
representation theory of the affine 
Weyl groups. 
For $k\geq 1$, they were
studied in \cite{Sta1, Athana1} among others.

Let $S:=S(V^{*} )$ be 
the symmetric algebra of the dual space $V^{*} $ of $V$.
Recall that the {\bf
cone}
\cite[Definition 1.15]{OT0}
$$
{\mathcal S}^{k}
:=
{\mathcal S}^{k} (\Phi^{+})
:=
{\mathbf c}{\mathit{Shi}}^{k}
$$
over 
$\mathit{Shi}^{k}$ 
is a central arrangement in 
an
$(\ell+1)$-dimensional Euclidean space $E:=\R^{\ell+1}$.
Choose $\alpha_{H} \in E^{*} $ with 
$H = \ker (\alpha_{H} )\in {\cal S}^{k} $. 
Let ${\Der} (S)$ denote the $S$-module of derivations of $S$
to itself.   
Since $V$  is embedded in $E$  as the affine hyperplane
defined by $z=1$,
we may consider the {\bf Ziegler restriction map} \cite{Z} 
\begin{align*}
{\mathbf{res}} :
D_{0} ({\mathcal S}^{k})
\longrightarrow
D(\A(\Phi), 2k)
\end{align*} 
by setting $z=0$.
Here 
\begin{align*}
&D_0({\mathcal S}^{k})
:=
\{
\theta\in \Der(S(E^{*} ))~|~
\theta(\alpha_H) \in \alpha_H  S(E^{*}) 
{\text{~for~each~}} H \in \mathcal{S}^{k},
\theta(z)=0 \},
\\
&D(\A(\Phi), 2k)
:=
\{\theta\in \Der(S(V^{*})) 
~|~ \theta(\alpha)\in\alpha^{2k}  S(V^{*} )
\text{~for each ~}\alpha\in\Phi^{+}
\}.
\end{align*} 
Yoshinaga \cite{Y}
verified the 
Edelman-Reiner conjecture
in \cite{EdRei}
by proving the surjectivity
of the Ziegler restriction map.
In particular,
the homogeneous part 
$$
{\mathbf{res}}
:
D_0({\mathcal S}^{k})_{kh} 
\longrightarrow
D(\A(\Phi), 2k)_{kh}
$$ 
of degree $kh$ of the Ziegler restriction map
is a linear isomorphism, where
$h$ denotes the Coxeter number.
On the other hand,
as we will
see in Proposition \ref{Xi},
$D(\A(\Phi), 2k)_{kh}$
 is isomorphic to
$V$ 
as
a  
$W$-module.

Now we are ready to introduce
the simple-root bases for 
$D_{0} ({\mathcal S}^{k})_{kh}$.

\begin{definition}
\label{SRB}
Pick a $W$-isomorphism
$
\Xi :
V \tilde{\longrightarrow}
D(\A(\Phi), 2k)_{kh}.
$
Define
a linear isomorphism
$\Theta:=  {\bf res}^{-1} \circ \Xi$:
\[
\xymatrix{
V  
\ar[dr]^(0.55){{\!\!\!{\scalebox{1.05}{$\Xi$}}}}_(0.55){\rotatebox{-30}{\scalebox{1.1}{\!\!\!\!\!\!$\sim$}}}
\ar[r]^(0.4){\scalebox{1.05}{$\Theta$}}_(0.4){\scalebox{1.05}{$\sim$}}
& 
D_{0}({\Shi}^{k})_{kh} \ar[d]^{\scalebox{1.05}{${\mathbf{res}}$}}_{\rotatebox{90}{$\sim$}} 
\ar@{}[ld]|(0.3){\scalebox{1.05}{$~\circlearrowright$}}
\\
& 
D({\mathcal{A}}(\Phi),2k)_{kh}. \\
}
\]
Let $\Delta:=\{\alpha_{1}, \dots , \alpha_{\ell} \}\subset V^{*} $ 
be the set of simple
 roots.
Let $\Delta^{*} :=\{\alpha_{1}^{*} , \dots , \alpha_{\ell}^{*}  \}\subset V$ 
denote  the dual basis of $\Delta$. 
Then

(1)  
the derivations
$$
\varphi^{+}_{j} := \Theta(\alpha_{j}^{*} ) \,\,(j=1, \dots, \ell)
$$  
are called a {\bf simple-root basis plus (SRB$+$)}, 
and

(2) the derivations
$$
\varphi^{-}_{j} := \sum_{p=1}^{\ell}  I^{*}(\alpha_{j}, \alpha_{p}) 
\varphi^{+}_{p}\,\,(j=1, \dots, \ell)
$$  
are called a {\bf simple-root basis minus (SRB$-$)}.
Here $I^{*} $ is the natural inner product on $V^{*} $ 
induced from the inner product $I$ on $V$. 
\end{definition}

These two bases
are uniquely determined up to a nonzero constant multiple
because of Schur's lemma.
The bases have the following nice
properties:

\begin{theorem}
\label{characterization} 
(1)
Let
 $
\varphi^{+}_{1},
\dots,
\varphi^{+}_{\ell}
$ be an SRB$+$.  Then
each  $
\varphi^{+}_{i}(\alpha_{j}+kz)$ is divisible by $\alpha_{j}+kz$ 
whenever $i\neq j$.
Conversely, if derivations  $
\phi^{+}_{1},
\dots,
\phi^{+}_{\ell}
\in
D_{0}({\Shi}^{k})_{kh}
$ are given and
each  $\phi^{+}_{i}(\alpha_{j}+kz)$ is divisible by $\alpha_{j}+kz$ 
whenever $i\neq j$, 
then
there exist nonzero constants
$
c^{+}_{1},
\dots,
c^{+}_{\ell}
$ 
satisfying
$\phi^{+}_{i} = c^{+}_{i} \varphi^{+}_{i}$ for any $i$. 

\noindent
(2) Let $
\varphi^{-}_{1},
\dots,
\varphi^{-}_{\ell}
$ be an
SRB$-$.
 Then
each  derivation $\varphi^{-}_{i}$ is divisible by $\alpha_{i}-kz$.
Conversely, if derivations  
$
\phi^{-}_{1},
\dots,
\phi^{-}_{\ell}
\in
D_{0}({\Shi}^{k})_{kh}
$ are given and
each derivation $\phi^{-}_{i}$ is divisible by $\alpha_{i}-kz$, 
then
there exist nonzero constants
$
c^{-}_{1},
\dots,
c^{-}_{\ell}
$ 
satisfying
$\phi^{-}_{i} = c^{-}_{i} \varphi^{-}_{i}$ for any $i$. 
\end{theorem}

The organization of this article is as follows. 
In Section 2, we review a recent refinement \cite{AY2} 
of Yoshinaga's freeness criterion \cite{Y} 
before proving the two key results (Propositions \ref{BGamma+}
and \ref{BGamma-})
which we apply in  Section 3 when
we prove  
Theorem \ref{characterization}.
We also characterize the simple roots
in terms of the freeness of deleted/added 
Shi arrangements in Theorem \ref{simplefree}.
The actions of a simple reflection on the SRB$+/-$ 
are studied as well.
In Section 4, we will describe a unique 
$W$-invariant derivation 
({\bf $k$-Euler derivation})
related to the Catalan arrangement in terms of the
SRB$+$.

\medskip

\noindent
\textbf{Acknowledgements}. 
The first author is partially supported by JSPS Grants-in-Aid for Young Scientists
(B)
No. 24740012. 
The second author is partially supported by JSPS Grants-in-Aid, Scientific Research
(A) 
No. 24244001.

\section{Freeness criteria}
In the rest of the article we use \cite{OT0} as a general reference.

Let $\cal C$ be a central arrangement in an $(\ell+1)$-dimensional vector space
$E$.
Choose $\alpha_{H} \in E^{*} $
with $\ker \alpha_{H} = H\in\cal C$.
Let $\mathbf m$ be a {\bf
multiplicity} of $\cal C$:
\[
{\mathbf m} : {\cal C} \rightarrow {\Z_{>0}}. 
\]
Let $S(E^{*})$ be the ring of polynomial functions on $E$.
Let $\Der(S(E^{*}))$ be the set of derivations of $S(E^{*})$ to itself:
\begin{align*} 
\Der(S(E^{*})):= \{\theta : S(E^{*})\rightarrow 
S(E^{*})~|~&\theta \text{~is $\R$-linear
and~} \theta(fg)=f\theta(g)+g\theta(f)\\
 &\text{~for~any~} f, g\in S(E^{*}) \}.
\end{align*} 
A derivation $\theta\in \Der(S(E^{*}))$
is said to be {\bf
homogeneous of degree} $d$ if
$\theta(\alpha)$ is a 
homogeneous polynomial 
of 
degree $d$
for any $\alpha\in V^{*} $  unless $\theta(\alpha)=0.$ 
Define an $S(E^{*})$-module
\[
D({\cal C}, {\mathbf m})
:=
\{
\theta\in\Der(S(E^{*}))
~|~
\theta(\alpha_{H} )\in \alpha_{H}^{{\bf m}(H)} S(E^{*})
\text{~for each~}H\in {\cal C}   
\}.
\]
When 
$D({\cal C}, {\mathbf m})
$   
is a free 
$S(E^{*} )$-module, we say that a multiarrangement
$({\cal C}, {\bf m}  )$ is {\bf free}.
We say that ${\cal C}$ is a {\bf free arrangement}
if $({\cal C}, {\bf 1}  )$ is free.
Here ${\bf 1} $ indicates the constant multiplicity whose value
is equal to one. 
When 
$({\cal C}, {\bf m}  )$ is {free},
$\mathbf{\exp}
({\cal C}, {\bf m}  )
$
of {\bf exponents} denotes the set of degrees of homogeneous basis for
$D({\cal C}, {\mathbf m})
$.
We simply write 
   ${\exp}({\cal C})$
instead of 
$\exp
({\cal C}, {\bf 1}  )
$
if $\mathcal C$ is a free arrangement.

For a fixed
hyperplane $H_{0}\in\calC $,
define a multiarrangement
$(\calC'', \bfz)$, which we call the {\bf Ziegler restriction}
\cite{Z}, by
\begin{align*} 
\calC'' :=
\{
H_{0} \cap K ~|~ K\in \calC':=\calC\setminus\{H_{0} \}\},
\,\,
\bfz (X):= |\{
K\in\calC'~|~ X=K\cap H_{0} 
\}|,
\end{align*} 
where $\calC''$ is an arrangement  in $H_{0} $ 
and
$X\in\calC''$. 
For the intersection lattice $L(\calC)$
\cite[Definition 2.1]{OT0}  of $\calC$ 
and any $Y\in L(\calC)$ 
define the {\bf
localization} $\calC_{Y} $ of $\calC$ at $Y$ by
$
\calC_{Y} :=\{H\in\calC~|~Y\subseteq H\}.
$
Let us present a recent refinement
of Yoshinaga's freeness criterion in
\cite{Y}:

\begin{theorem}
[\cite{AY2}]
\label{AbeYoshinagaCriterion}
Suppose $\ell+1>3$.
 For a central arrangement $\calC$ and an 
arbitrary hyperplane
$H_{0}\in\calC $, the following two conditions are
equivalent:

(1) $\calC$ is a free arrangement,

(2) (2-i) the Ziegler restriction $(\calC'', \bfz)$  
is free
and
(2-ii) 
$\calC_{Y}$ is free for any $Y\in L(\calC)$
 such that $Y\subset H_{0} $ with $\codim_{H_{0}} Y =2$. 
\end{theorem}

For a fixed hyperplane $H_{0} \in\mathcal C$, we may choose 
a basis $x_{1}, x_{2},\dots , x_{\ell}, z  $ 
for $E^{*} $ so that the hyperplane $H_{0} $ is defined by
the equation $z=0$.
Then the Ziegler restriction 
$(\mathcal C'', \mathbf z)$ is a multiarrangement in $H_{0} $.
Let
\[
D_{0}(\mathcal C):=\{\theta\in D(\mathcal C)~|~
\theta(z)=0
\}. 
\]
Then
\[
D(\mathcal C)
=
S\theta_{E} \oplus 
D_{0} (\mathcal C),
\]
where
$\theta_{E}$ is the Euler derivation.
Note that $D_{0} (\mathcal C)$ is a free $S(E^{*} )$-module
if and only if
$\mathcal C$ is a free arrangement. 
When 
$\mathcal C$ is a free arrangement,
let $\exp_{0} (\mathcal C)$ 
denote the set of degrees of homogeneous basis for
$D_{0}(\mathcal C) $. 
Note that the set  $\exp_{0} (\mathcal C)$ 
does not depend upon the choice of $H_{0}. $ 
When we describe $\exp_{0} (\mathcal C) $,
 we will
use the notation $a^{n} $ instead of listing $a, \dots, a$ ($n$ times).

\begin{theorem}
[Ziegler \cite{Z}]
\label{Ziegler}
The Ziegler restriction  map
\[
{\mathbf{res}}:
D_{0} (\mathcal C)
\rightarrow
D(\calC'', \bfz)
\]
defined by setting $z=0$
is surjective if   
$\mathcal C$ is
a free arrangement.
\end{theorem}

The following theorem was proved 
in
\cite{AY} using
the shift isomorphism of 
Coxeter multiarrangements:

\begin{theorem}[\cite{AY}, Corollary 12]
\label{AbeYoshinagaCorollary}
Let $\A$ be a Coxeter arrangement in an $\ell$-dimensional Euclidean space. 
For a $\{0, 1\}$-valued multiplicity $\mathbf{m} :
\A \rightarrow \{0, 1\}$ and an integer $k>0$, the following conditions are
equivalent:

(1) a multiarrangement $(\A, \mathbf{m} )$ is free with 
exponents $(e_{1} , \dots, e_{\ell} )$.

(2) a multiarrangement $(\A, 2k+\mathbf{m} )$ is free with 
exponents $(kh+e_{1} , \dots, kh+e_{\ell} )$.

(3) a multiarrangement $(\A, 2k-\mathbf{m} )$ is free with 
exponents $(kh-e_{1} , \dots, kh-e_{\ell} )$.
\end{theorem}

Now we go back to the situation in Section 1:
let $\Phi, \Phi^{+}$ and $ \Delta $ be an irreducible root system,
a set of postive roots, and the set of simple roots
respectively.
The following two Propositions \ref{BGamma+} and
\ref{BGamma-} are keys to our proof of
Theorem \ref{characterization}.   They are dual to each other. 
 
\begin{prop}
\label{BGamma+} 
For any subset $\Gamma$ of the simple system
$\Delta$, the arrangement
$$
\B_{\Gamma}^{+} :=
\B_{\Gamma}^{+}(\Phi^{+}) :=
{\mathcal S}^{k} 
\cup \bigcup_{\alpha\in \Gamma} \{{\mathbf c} H_{\alpha, -k }\}
$$
is a free arrangement
with 
$${\rm exp}_{0} (\B_{\Gamma}^{+})
=
((kh+1)^{|\Gamma|}, (kh)^{\ell-|\Gamma|}).
$$ 
\end{prop}

\noindent
\textbf{proof}.
{\it Case 1.}
When $\ell=2$, 
$\Phi$ is of the type either $A_{2} $, $B_{2} $ 
or $G_{2} $.
Then $\exp_{0} (\mathcal S^{k})=((kh)^{2})
=(kh,kh)$ and 
$\Delta=\{\alpha_1,\alpha_2\}$.
For an affine $2$-arrangement $\A$ and an affine line $H_{0}$,
define
\[
\A\cap H_{0} :=
\{K\cap H_{0} \mid K\in \A, K\neq H_{0} \}.
\]  
 Then, by directly counting intersection points, 
we get the following equalities:
\begin{align*} 
|\mathit{Shi}^{k}  
\cap H_{\alpha, -k}|&=kh~~(\alpha\in\Delta),\\
|(\mathit{Shi}^{k} \cup \{H_{\alpha_1,-k}\}) 
\cap H_{\alpha_{2} , -k}|&=kh+1.
\end{align*} 
Thus
we 
may verify the statement by
applying the addition theorem \cite{T} \cite[Theorem 4.49]{OT0}
to $\mathcal S^{k} $ 
for the types of
$A_{2}, B_{2}$ and $G_{2}$. 
\medskip

{\it Case 2.}
Suppose that $\ell \ge 3$. 
We will apply Theorem \ref{AbeYoshinagaCriterion}
by verifying the two conditions
(2-i) and (2-ii).

(2-i) 
Note that
the Ziegler restriction of $\B_\Gamma^{+}$ to
the hyperplane $H_\infty$
at infinity
coincides with $(\A(\Phi), {\bfz}_\Gamma^{+} )$,
where
\[
\bfz_{\Gamma}^{+} (H_{\alpha}) =
 |\{j ~|~ {\mathbf c}H_{\alpha, j} \in 
{\mathcal B}_{\Gamma}^{+}\}|
=
\begin{cases} 
2k+1 & {\text{~if~}} \alpha\in \Gamma\\
2k & {\text{~otherwise}}
\end{cases} 
\,\,\,
\,\,\,
\,\,\,
(\alpha\in\Phi^{+}).
\]
If $\Gamma$ is empty, then
$\bfz_{\Gamma}^{+}  = \bfz_{\emptyset}^{+}  \equiv 2k$. 
Note that $\Gamma$ is linearly independent
because it is a set consisting of simple roots.
Thus the arrangement
\[
\A(\Gamma) := \{H_{\alpha} ~|~ \alpha\in \Gamma\}
\]
is a free (Boolean)
subarrangement of $\A(\Phi)$. 
Let $\chi_{\Gamma} $ be
the characteristic function of $\A(\Gamma)$ in $\A(\Phi)$:
\[
\chi_{\Gamma} 
(H_{\alpha} ) =
\begin{cases} 
1 & {\text{~if~}} \alpha \in \Gamma\\
0 & {\text{~otherwise}}
\end{cases} 
\,\,\,
\,\,\,
\,\,\,
(\alpha\in\Phi^{+}).
\]
 Since $\bfz_{\Gamma}^{+}  =
\bfz_{\emptyset}^{+}  + \chi_{\Gamma} $, we may apply
Theorem \ref{AbeYoshinagaCorollary} 
 to conclude
 that 
$(\A(\Phi), \bfz_\Gamma)$ 
is a free multiarrangement
with exponents $((kh+1)^{|\Gamma|},(kh)^{\ell-|\Gamma|})$.

(2-ii) 
We will prove that 
$(\B_\Gamma^{+})_{Y}$ is free for any
$Y \in L(\B_\Gamma^{+})$
such that $Y \subset H_\infty$
with $\codim_{H_{\infty}} Y=2$. 
Define $X$ to be the unique subspace of $E$ 
such that 
 $X\in L(\A(\Phi))$ and
${\mathbf c}X \cap H_{\infty} = Y$. 
Let
$X^{\perp} :=
\{\alpha\in V^{*} ~|~ \alpha|_{X} \equiv 0\}$. 
Then
$
\Phi_{X} 
:=
\Phi\cap X^{\perp} 
$
 is also a (not necessarily irreducible)
root system in $X^{\perp} $.
The set
of positive roots 
of
$\Phi_{X} 
$ 
is induced from $\Phi^{+} $:
$\Phi_{X}^{+}=\Phi^{+} \cap \Phi_{X}$.
It is not hard to see (e.g. \cite[Lemma 3.1]{AT})
that
$$
({\mathcal S}^{k})_{Y} 
=
{\mathbf c}
\left(
\mathit{Shi}^{k}(\Phi^{+}_{X}) 
\times \emptyset_{X} \right),$$
where $\emptyset_{X} $ denotes the empty arrangement
in $X$. 
Since $\dim X^{\perp} = 2$,
$\Phi_{X}$ 
is either
of the type 
$A_{1} \times A_{1}$,
$A_{2} $, 
$B_{2} $ or $G_{2} $.

{\it Case 2.1.}
When 
$\Phi_{X}$ 
is
of the type 
$A_{1} \times A_{1}$,
the arrangement 
$
\mathit{Shi}^{k}(\Phi^{+}_{X})
$ 
is a product of two affine
$1$-arrangements.
Thus any subarrangement
of
$
({\mathcal S}^{k})_{Y} 
$
is a free arrangement.
In particular,
$(\B_{\Gamma}^{+})_{Y}  $ 
is a free arrangement.

{\it Case 2.2.}
Suppose
that
$\Phi_{X}$ 
is
of the type either
$A_{2}$,
$B_{2}$
or
$G_{2}$.
Suppose that $\alpha\in\Phi_{X} $ is a simple root of $\Phi$.
Then
$\alpha$
  is also a simple root of 
$\Phi_{X}$ because it cannot be expressed as a sum of two positive roots of
$\Phi_{X}$.
Thus $\Phi_{X}\cap\Gamma$ consists of simple roots of $\Phi_{X} $.
Therefore 
\begin{align*} 
(\B_{\Gamma}^{+})_{Y}
&
=
({\mathcal S}^{k})_{Y} 
\cup
\bigcup_{\alpha\in \Phi_{X} \cap \Gamma} \{\mathbf{c}H_{\alpha, k} \}\\
&
= 
{\mathbf c}
\left(
\left(
\mathit{Shi}^{k}(\Phi^{+}_{X}) 
\cup
\bigcup_{\alpha\in \Phi_{X} \cap \Gamma} \{H_{\alpha, k} \}
\right)
\times \emptyset_{X}
\right)
\end{align*} 
is a free arrangement
because of {\it Case 1}.

Now we apply 
Theorem \ref{AbeYoshinagaCriterion}
to complete the proof. 
\owari

\begin{cor}
The vector space $D_{0}(\B_{\Gamma}^{+})_{kh}  $ 
is $(\ell-|\Gamma|)$-dimensional. 
\label{BGammadim+}
\end{cor}

\begin{prop}
\label{BGamma-} 
For any subset $\Gamma$ of the simple system
$\Delta$, the arrangement
$$
\B_{\Gamma}^{-} :=
\B_{\Gamma}^{-}(\Phi^{+}) :=
{\mathcal S}^{k}
\setminus \bigcup_{\alpha\in \Gamma} \{{\mathbf c} H_{\alpha, k }\}
$$
is a free arrangement
with 
$${\rm exp}_{0} (\B_{\Gamma}^{-})
=
((kh-1)^{|\Gamma|}, (kh)^{\ell-|\Gamma|}).
$$ 
\end{prop}

\noindent
\textbf{proof}.
{\it Case 1.}
When $\ell=2$, 
$\Phi$ is of the type either $A_{2} $, $B_{2} $ 
or $G_{2} $.
Let $\exp_{0} (\mathcal S^{k})=((kh)^{2})$ and 
$\Delta=\{\alpha_1,\alpha_2\}$.
 Then, by directly counting intersection points, 
we get the following equalities:
\begin{align*} 
|\mathit{Shi}^{k}  
\cap H_{\alpha,k}|&=kh~~(\alpha\in\Delta),\\
|(\mathit{Shi}^{k} \setminus \{H_{\alpha_1,k}\}) 
\cap H_{\alpha_{2} , k}|&=kh-1.
\end{align*} 
Thus
we 
may verify the statement by
applying the deletion theorem \cite{T} \cite[Theorem 4.49]{OT0}
to $\mathcal S^{k} $ 
for the types of
$A_{2}, B_{2}$ and $G_{2}$. 

The rest is exactly the
same as the proof of
Proposition \ref{BGamma+} if one replaces
$\B_{\Gamma}^{+}$, 
$kh+1$, 
$H_{\bullet, -k}$, 
$2k+1$, 
$\cup$, 
${\mathbf z}_{\Gamma}^{+}  $, 
${{\mathbf z}}_{\Gamma}^{+}+\chi_{\Gamma}$ 
with
$\B_{\Gamma}^{-}$, 
$kh-1$, 
$H_{\bullet, k}$, 
$2k-1$, 
$\setminus$, 
${\mathbf z}_{\Gamma}^{-}  $, 
${\mathbf z}_{\Gamma}^{-}-\chi_{\Gamma}$
respectively. 
\owari

\begin{cor}
The vector space $D_{0}(\B_{\Gamma}^{-})_{kh-1}  $ 
is
$|\Gamma|$-dimensional.
\label{BGammadim-}
\end{cor}

\section{Proof of main results} 

We will prove Theorem \ref{characterization} in this section.
Fix $k\in\Z_{>0} $ throughout in the rest of this article.
We first see that $V$ and $D(\A(\Phi), 2k)_{kh} $ are 
$W$-isomorphic as mentioned in Section 1.
Let $F$ be the field of quotients of $S=S(V^{*})
=\R[x_{1} , \dots , x_{\ell}]$.
Recall a primitive derivation $D\in \Der(F)$ associated with $\A(\Phi)$:
$D$ satisfies
\[
D(P_{j} ) = 
\begin{cases} 
c \in \R^{\times} &\text{~if $j=\ell$ },\\ 
0                & \text{~if $1\leq j\leq \ell-1$ }.
\end{cases} 
\]
Here $P_{1}, \dots , P_{\ell}  $ 
are basic invariants of the invariant subring
$S^{W}  $ with
\[
2=\deg P_{1} < \deg P_{2} \leq \dots\leq \deg P_{\ell-1} < \deg P_{\ell} =h.
\]
Then the choice of $D$ has the ambiguity of nonzero constant multiples. 

Consider the Levi-Civita connection with respect to the inner product
$I$:
\begin{align*} 
\nabla : 
{\rm Der}(F) \times {\rm Der}(F)  \longrightarrow
{\rm Der}(F),
\,\,\,\,
(\xi, \eta)  \mapsto
\nabla_{\xi} \eta
\end{align*} 
Note that
$$
\left(\nabla_{\xi}\eta\right)(\alpha) 
= 
\xi(\eta(\alpha))
\,\,\,(\xi, \eta\in \Der(F), \alpha\in V^{*} )
$$
because $I: V\times V \rightarrow \R$ is real
number-valued.
Consider $T:= \R[P_{1} ,\dots, P_{\ell-1}]$-linear
covariant derivative
$
\nabla_{D} : \Der(F) \rightarrow \Der(F).
$
By
\cite{AT10} it induces a $T$-linear bijection 
\[
\nabla_{D} 
:
D(\A(\Phi), 2k+1)^{W} 
\tilde{\longrightarrow}
D(\A(\Phi), 2k-1)^{W}~~(k>0).
\]
The covariant derivative $\nabla_{D} $ was introduced by K. Saito (e.g. \cite{Sa93}) to study the flat structure (or the Frobenius manifold structure) of the
orbit space $V/W$.
Let $\theta_{E}:=\sum_{i=1}^{\ell} x_{i} (\partial/\partial x_{i}) 
$ denote the Euler derivation.  Since $\theta_{E}\in D(\A(\Phi), 1)^{W}$,
one has
\[
\nabla_{D}^{-k} \theta_{E} \in D(\A(\Phi), 2k+1)^{W}  
\]
 which plays a principal role in this section.

For any $v\in V$, there exists a unique derivation 
$\partial_{v} \in \Der(S)_{0} $ of degree zero
such that $$
\partial_{v} (\alpha) := \left<\alpha, v
\right>\,\,\,\,(\alpha\in V^{*} ).
$$  
Thus we may identify
$V$ with $\Der(S)_{0} $ by the $W$-isomorphism
\[
V \longrightarrow \Der(S)_{0}
\]
defined by
$
v\mapsto \partial_{v}.  
$  

Let $\Omega(S)$ be the $S$-module of regular one-forms:
\[
\Omega(S)=S(dx_{1} )\oplus\dots\oplus S(dx_{\ell}).
\]
 We may identify $V^{*} $ with $\Omega(S)_{0} $ 
as $W$-modules by the bijection $\alpha \mapsto d\alpha$.  

Recall the $W$-invariant dual inner product 
$
I^{*} : V^{*} \times V^{*} \rightarrow \R.
$ 
Define a $W$-isomorphism
\[
I^{*} : \Omega(S)_{0}  \rightarrow \Der(S)_{0} 
\]
 by
$
\left(I^{*} (d\alpha)\right) 
(\beta):= I^{*} (d\alpha, d\beta) 
\,\,\,\,(\alpha, \beta\in V^{*}).
$ 

When a $W$-isomorphism $\Xi: \Der(S)_{0} 
 \rightarrow D(\A(\Phi), 2k)_{kh} $ is given,
we have the following commutative diagram in which the new maps
$\Xi^{*} $
and
$\Theta^{*} $ are defined:

$$
\xymatrix
{
&
{
D_{0}({\Shi}^{k})_{kh}}
\ar[ddd]^{\scalebox{1.1}{${\bf{res}}$}}_{\rotatebox{270}{\scalebox{3.0}[1.0]{$\widetilde{}$}}}
&
\\
&
{\Theta^{*}}~~~~~~~~~~~~
~~~~~~~~{\Theta}
&
\\
&
&
\\
&
{ D(\mathcal{A}(\Phi), 2k)_{kh}}
&
\\
&
\rotatebox{180}{\scalebox{1.3}{$\circlearrowright$}}
&
\\
{V^{*}=\Omega(S)_{0}
\ar[rr]_{\scalebox{1.1}{$I^{*}$}}}
\ar[ruu]^{\scalebox{1.1}{$\Xi^{*}$}}_{
\rotatebox{30}{\scalebox{3.0}[1.0]
{$\widetilde{}$}}}
\ar[ruuuuu]^{\scalebox{1.1} {}}_(0.53){{\rotatebox{180}{{\scalebox{1.3}{$\circlearrowright$}}}}}
\ar@{}[ru]|{\scalebox{1.1}{$W$-iso.}} 
&
&
{ 
{\rm Der}(S)_{0}=V}
\ar[luu]^{\rotatebox{330}{\scalebox{3.0}[1.0]{$\widetilde{}$}}}_{
\scalebox{1.1}{$\Xi$}} 
\ar[luuuuu]_{\scalebox{1.1} {}}^(0.53){\rotatebox{180}{\scalebox{1.3}{$\circlearrowright$}}}
\ar@{}[lu]|{\scalebox{1.1}{$W$-iso.}}
}
$$

\begin{prop}
\label{nablaD} 
(1) For any primitive derivation $D$, define
\label{Xi} 
\[
\Xi_{D}  : \Der(S)_{0}  {\longrightarrow} D(\A(\Phi), 2k)_{kh}  
\]
by
\[
\Xi_{D} (\partial_{v}) 
:= 
\nabla_{\partial_{v}} \nabla_{D}^{-k} \theta_{E}
\,\,\,\,(v\in V). 
\]
Then $\Xi_{D} $ is a $W$-isomorphism. 

\noindent
(2) Conversely, for any $W$-isomorphism
$
\Xi  : \Der(S)_{0}  \tilde{\longrightarrow} D(\A(\Phi), 2k)_{kh}  
$,
there exists a unique primitive derivation
such that $\Xi = \Xi_{D} $. 
 \end{prop} 

\noindent
{\bf proof}.
(1) was proved by Yoshinaga in \cite{Y02}.
(See \cite{T02} also.)
(2) follows from (1) and Schur's lemma. 
\owari 

\begin{remark}
\label{remark1.3}

(1)
The  linear isomorphism $\Theta$ 
has the ambiguity of nonzero constant multiples
as in the case of the choice of
$\Xi $.
By Proposition \ref{nablaD}, the map $\Theta$ is uniquely determined 
when a primitive derivation 
$D$ is fixed.

\noindent
(2)
Note that
the arrangement
$
{\Shi}^{k}
$  
is not $W$-stable.
Therefore 
the $\ell$-dimensional vector space  $D_{0}({\Shi}^{k})_{kh}$
is not naturally a $W$-module, while the $\ell$-dimensional vector spaces
$V$ and $D({\mathcal{A}}(\Phi),2k)_{kh}$ are both
$W$-modules. 
\end{remark}


\bigskip

\noindent
{\bf{proof of Theorem \ref{characterization}.}}
Let $\Delta := \{\alpha_{1}, \dots, \alpha_{\ell}  \}$
be the set of simple roots.  
Let $\{\alpha_{1}^{*} , \dots, \alpha_{\ell}^{*}   \}$
be the dual basis to $\Delta$: $\left<
\alpha_{i}, \alpha^{*}_{j}   \right> = \delta_{ij} $ 
(Kronecker's delta).
Fix a primitive derivation $D$ 
and let $\Xi = \Xi_{D} $.   

(1) Let $\varphi_{i}^{+} := \Theta (\partial_{\alpha_{i}^{*}  } )
\,\,\,
(1\leq i\leq \ell) 
$
be the SRB$+$. 
We have
\begin{align*} 
\left[{\bf{res}}(\varphi_{i}^{+})\right] (\alpha_{j})
&=
\left[{\bf{res}}(\Theta(\partial_{\alpha_{i}^{*} 
}))\right] 
(\alpha_{j})
=
\left[\Xi(\partial_{\alpha_{i}^{*} })\right] 
(\alpha_{j})
=
(
\nabla_{\partial_{\alpha_{i}^{*} }}\nabla_{D}^{-k} \theta_{E}  
)
(\alpha_{j} )\\
&=
\partial_{\alpha_{i}^{*} }\left(\left(\nabla_{D}^{-k} \theta_{E}
\right)
(\alpha_{j} )\right)
\in
\alpha_{j}^{2k+1} S
 \end{align*} 
because $\partial_{{\alpha}_{i}^{*}}(\alpha_{j})=0$ if $j\neq i$.
Define $\Gamma^{+}_{i}:=\Delta\setminus\{\alpha_{i} \}$. 
Then we have the following commutative diagram
\[
\xymatrix{ 
D_{0} ({\mathcal S}^{k})_{kh} 
\ar[r]^(0.43){\bf {res}}_(0.43){\scalebox{1.5}{$\widetilde{}$}}  
& 
D (\mathcal{A}(\Phi),2k)_{kh}
\\
D_{0} ({\mathcal B}_{\Gamma^{+}_{i}  }^{+})_{kh} 
\ar[r]^(0.43){\bf {res}}_(0.43){\scalebox{1.5}{$\widetilde{}$}}  
\ar@{}[u]|{\scalebox{1.3}{$\bigcup$}}
& 
D(\mathcal{A}(\Phi),{\bf z}^{+}_{\Gamma^{+}_{i}})_{kh} 
\ar@{}[u]|{\scalebox{1.3}{$\bigcup$}}
}
\]
because of Proposition \ref{BGamma+}
and Theorem \ref{Ziegler}. 
Since 
$\varphi_{i}^{+}\in  D_{0} (S^{k})_{kh} $ 
and 
${\bf{res}}\left(\varphi_{i}^{+}\right)
\in  
D(\mathcal{A}(\Phi),{\bf z}^{+}_{\Gamma^{+}_{i}})_{kh}
$,
we conclude  that
$\varphi_{i}^{+}\in  
D_{0} ({\mathcal B}_{\Gamma^{+}_{i}  }^{+})_{kh} 
$
by chasing the diagram above. 

The ``uniqueness part'' 
follows from the equality
$$\dim D_{0} ({\mathcal B}_{\Gamma^{+}_{i}  }^{+})_{kh} =1$$ 
which is a consequence of
Corollary \ref{BGammadim+} because 
$|\Gamma^{+}_{i}|=\ell-1$.

(2) Recall 
the SRB$-$:
\begin{align*} 
\varphi_{i}^{-} 
&:= \sum_{p=1}^{\ell}  
I^{*}(\alpha_{i}, \alpha_{p}) \varphi_{p}^{+}
=
\Theta (
\sum_{p=1}^{\ell}  
I^{*}(d\alpha_{i}, d\alpha_{p}) \partial_{\alpha_{p}^{*} }
)
=
\Theta (
I^{*}(d\alpha_{i})
)
\,\,\,\,
(1\leq i\leq \ell).
\end{align*} 
We have
\begin{align*} 
{\bf{res}}(\varphi_{i}^{-}) 
&=
{\bf{res}}\left(\Theta (I^{*}(d\alpha_{i}))\right)
=
\Xi(I^{*}(d\alpha_{i}))
=
\nabla_{I^{*}(d\alpha_{i})}\nabla_{D}^{-k} \theta_{E}.
 \end{align*} 
Let $s_{i} $ denote the orthogonal reflection
with respect to $\alpha_{i} $. 
Since
$s_{i} (d\alpha_{i} )= -d\alpha_{i} $,
we have
\[
s_{i} (
{\bf{res}}(\varphi_{i}^{-}) 
)
=
-
{\bf{res}}(\varphi_{i}^{-}). 
\] 
Express
\[
{\bf{res}}(\varphi_{i}^{-}) = \sum_{p=1}^{\ell} f_{p} 
\partial_{\alpha_{p}^{*}  }
\,\,\,\,\,(f_{p} \in S).  
\]
Recall $s_{i} (\partial_{\alpha^{*}_{p} }) = \partial_{\alpha^{*}_{p} }$ 
whenever $p\neq i$. 
Thus we have 
$s_{i} (f_{p})= -f_{p} $ 
whenever $p\neq i$. 
Therefore 
$f_{p}$ is divisible by $\alpha_{i}$ 
whenever $p\neq i$. 
We also know that
\[
f_{i}
= 
\left[{\bf{res}}(\varphi_{i}^{-})\right]
(\alpha_{i} ) 
=
(I^{*}(d\alpha_{i}))
\left(
(\nabla_{D}^{-k} \theta_{E})
(\alpha_{i} )\right) 
\]
is divisible by $\alpha_{i}^{2k}$
because $\nabla_{D}^{-k} \theta_{E} (\alpha_{i} ) $ is divisible by
$\alpha_{i}^{2k+1}$. 
Therefore 
we conclude that $f_{i} $ is divisible by 
$\alpha_{i} $ for any $i$. 
Define $\Gamma^{-}_{i}:=\{\alpha_{i} \}$. 
Then we have the following commutative diagram
\[
\xymatrix{ 
D_{0} ({\mathcal S}^{k})_{kh} 
\ar[r]^{\bf {res}}_{\scalebox{1.5}{$\widetilde{}$}}  
& 
D (\mathcal{A}(\Phi),2k)_{kh}
\\
(\alpha_{i} -kz)\cdot
D_{0} ({\mathcal B}_{\Gamma^{-}_{i}  }^{-})_{kh-1} 
\ar[r]^(0.46){\bf {~~~~~~res}}_(0.46){\scalebox{1.5}
{~~~~$\widetilde{}$}}  
\ar@{}[u]|{\scalebox{1.3}{$\bigcup$}}
& 
\alpha_{i} \cdot D(\mathcal{A}(\Phi),{\bf z}^{-}_{\Gamma_{i}^{-}  })_{kh-1} 
\ar@{}[u]|{\scalebox{1.3}{$\bigcup$}}
}
\]
because of Proposition \ref{BGamma-}
and Theorem \ref{Ziegler}. 
Since 
$\varphi_{i}^{-}\in  D_{0} (S^{k})_{kh} $ 
and 
$$
{\bf{res}}\left(\varphi_{i}^{-}\right)
\in  
\alpha_{i} \cdot D(\mathcal{A}(\Phi),{\bf z}^{-}_{\Gamma^{-}_{i}})_{kh-1}
,$$
we may conclude  
$\varphi_{i}^{-}\in  
(\alpha_{i} - kz)\cdot
D_{0} ({\mathcal B}_{\Gamma^{-}_{i}}^{-})_{kh-1} 
$
by chasing the diagram above. 

The ``uniqueness part'' follows from 
the equality
$$\dim D_{0} ({\mathcal B}_{\Gamma^{-}_{i}  }^{-})_{kh-1} =1$$ 
which is a consequence from
Corollary \ref{BGammadim-} because 
$|\Gamma^{-}_{i}|=1$. 
\owari

\medskip

The following proposition asserts that the simple roots 
can be characterized by the freeness of an added/deleted 
Shi arrangement:

\begin{theorem}
Let $\alpha\in \Phi^{+}$. Then

\noindent
(1)
the arrangement
$ {\mathcal S}^{k} \cup \{\mathbf{c}H_{\alpha, -k} \}$ is a
free 
arrangement
if and only if $\alpha$ is a simple root,
and

\noindent
(2)
the
arrangement
$ {\mathcal S}^{k} \setminus \{\mathbf{c}H_{\alpha, k} \}$ is a
free arrangement
if and only if $\alpha$ is a simple root.
\label{simplefree}
\end{theorem}

\noindent
\textbf{proof}.
(1)
By Proposition \ref{BGamma+}, 
the ``if part'' is already proved.
Assume that $\alpha\in \Phi^{+}$ is a non-simple root.
We will prove that
$\calS^{k} \cup \{\mathbf{c}H_{\alpha, -k} \}$
is not free.
We may express
$\alpha=\beta_1+\beta_2$
with 
$\beta_{i} \in \Phi^{+} \,\,(i=1, 2)$.
Let $H_i$ be the 
hyperplane defined
by
$\beta_i = 0  \,\,(i=1, 2)$.
Then  $X:=H_1 \cap H_2$ 
is of codimension two in $V$ because of 
basic properties of the root systems.
As in the proof of Proposition \ref{BGamma+}, 
consider the two-dimensional
root system
$\Phi_X = \Phi \cap X^{\perp} $.
Note that 
$\Phi_{X} $ is not of the type
$A_{1} \times A_{1} $ because
$\{\alpha, \beta_{1} , \beta_{2} \}\subseteq \Phi_{X} $.  
Recall that the localization
$\left(\calS^{k} \cup \{\mathbf{c}H_{\alpha, -k} \}\right)_{Y}$ 
is free if 
$\calS^{k} \cup \{\mathbf{c}H_{\alpha, -k} \}$
is free.
In the root system $\Phi_{X} $, $\alpha$ cannot be a simple 
root since $\alpha$ is a sum of two positive roots
$\beta_{1} $ and $\beta_{2} $.
Thus we may assume that
$\Phi$ is a two-dimensional root system 
without loss of generality.
In this case 
the arrangement
$\calS^{k} \cup \{\mathbf{c}H_{\alpha, -k} \}$
is not free because of the addition theorem and
the equality
$$
|{Shi}^{k}  
\cap H_{\alpha, -k}|=kh+1~~(\alpha\not\in\Delta),
$$
which can be verified by directly
counting intersection points.

(2)
Exactly the
same as the proof of
(1)
 if one replaces
Proposition \ref{BGamma+},
$\cup$, 
$kh+1$, 
$H_{\bullet, -k}$, 
addition theorem
with
Proposition
 \ref{BGamma-},
$\setminus$, 
$kh-1$, 
$H_{\bullet, k}$, 
deletion theorem
respectively. 
\owari

\bigskip

Since the SRB bases uniquely exist,
it is 
natural to pose the following problem:

\begin{problem}
Give an explicit construction of
the SRB$+$ and SRB$-$
for
$D_{0} ({\mathcal S}^{k})$ 
for every irreducible root system and every 
$k\in \Z_{>0} $. 
\end{problem}

This problem has been solved in
\cite{SuT, Su, GPT} for the SRB$-$ 
when $k=1$ and
the root system is either of the type $A_{\ell},
B_{\ell}, C_{\ell}, D_{\ell}    $ or
$G_{2} $. The SRB$+$ can be obtained from the SRB$-$
by Definition \ref{SRB} (2). 

\medskip

In the rest of this section we will describe the action of a 
simple reflection on elements of the SRB$+/-$.

\begin{theorem}
\label{simplereflectionaction}
Let $\Delta = \{\alpha_{1}, \dots, \alpha_{\ell}  \}$ 
be the set of simple roots.
Let $s_{i} $ denote the simple reflection
with respect to a simple root $\alpha_{i} $ 
with $s_{i}(\alpha_{i} )=-\alpha_{i}  
\,\,(1\leq i\leq \ell)$.  

(1) If $1\leq j\leq \ell$ with $i\neq j$, then 
$s_{i} (\varphi_{j}^{+}  ) = \varphi_{j}^{+}  $.

(2) If we express 
$
\varphi_{i}^{-} 
=
(\alpha_{i} - kz) 
\hat{\varphi}_{i}^{-}, 
$
then
$s_{i}(\hat{\varphi}_{i}^{-}) 
=
\hat{\varphi}_{i}^{-}.
$
 \end{theorem}
 
\noindent
{\bf proof.}
Recall that 
(e.g., see \cite[Ch. VI, \S 1, no. 6, Corollary 1]{bou})
\[
s_{i} 
(
\Phi^{+} \setminus \{\alpha_{i} \}
)
=
\Phi^{+} \setminus \{\alpha_{i} \}.
\]

(1) When $i\neq j$, 
$\varphi^{+}_{j}(\alpha_{i} + kz)  $ 
is divisible by 
$\alpha_{i}+kz $ because of Theorem \ref{characterization} (1).
Thus
 $\left(
s_{i} (\varphi^{+}_{j})
\right)(\alpha_{i} - kz)  $ 
is divisible by 
$\alpha_{i}-kz $.
Using 
$
s_{i} 
(
\Phi^{+} \setminus \{\alpha_{i} \}
)
=
\Phi^{+} \setminus \{\alpha_{i} \},
$ 
we may conclude that
$
s_{i} (\varphi^{+}_{j})
\in
D_{0} (\calS^{k} ).
$
Consider the Ziegler restriction map
\[
\xymatrix{ 
D_{0} ({\mathcal S}^{k})_{kh} 
\ar[r]^(0.43){\bf {res}}_(0.43){\scalebox{1.5}{$\widetilde{}$}}  
& 
D (\mathcal{A}(\Phi),2k)_{kh}
}
\]
Then we obtain
\begin{align*} 
{\mathbf{res} } (s_{i}(\varphi_{j}^{+}  ) )
=
s_{i}({\mathbf{res}} (\varphi_{j}^{+})  )
=
s_{i}(\Xi(\partial_{\alpha_{j}^{*} })  )
=
\Xi(s_{i}(\partial_{\alpha_{j}^{*} })  )
=
\Xi(\partial_{\alpha_{j}^{*} })
=
{\mathbf{res} } (\varphi_{j}^{+})
\end{align*} 
because 
$
s_{i}(\partial_{\alpha_{j}^{*} })
=
\partial_{\alpha_{j}^{*} }
$ 
whenever $i\neq j$.
Hence we have $s_{i}(\varphi_{j}^{+}  )=
\varphi_{j}^{+}. $  

(2) Define $\Gamma_{i}^{-} := \{\alpha_{i} \}$ 
as in the proof of Theorem \ref{characterization} (2).
Then $\hat{\varphi}_{i}^{-}\in D_{0}
(\B^{-}_{\Gamma_{i}^{-}}    )_{kh-1}.$ 
Since
$$
\exp_{0} (\B^{-}_{\Gamma_{i}^{-}})
=
(
kh-1, \left(kh\right)^{\ell-1} 
),$$ 
the Ziegler restriction map
\[
\xymatrix{ 
D_{0} (\B^{-}_{\Gamma^{-}_{i}  })_{kh-1} 
\ar[r]^(0.43){
\!\!\!\!\!\!\!\bf {res}}_(0.43){\scalebox{1.5}{
\!\!\!\!\!\!\!$\widetilde{}$}}  
& 
D (\mathcal{A}(\Phi),2k-\chi_{\Gamma_{i}^{-}  } )_{kh-1}
}
\]
is bijective.
Note that the arrangement
$
\B^{-}_{\Gamma^{-}_{i}  }  
$ 
is $s_{i} $-stable because
\[
\B^{-}_{\Gamma^{-}_{i}  }  
=
\bigcup_{\beta\in\Phi^{+}\setminus \{\alpha_{i}\}}  
\{
H_{\beta, p } ~|~
1-k\leq p\leq k
\}
\cup
\{
H_{\alpha_{i}, p } ~|~
1-k\leq p\leq k-1
\}.
\]
 Then 
$$s_{i} (
\hat{\varphi}^{-}_{i}   
)
\in 
D_{0} 
(
\B^{-}_{\Gamma^{-}_{i}  }  
).
$$
We also verify
that
\[
{\mathbf{res}} (
\hat{\varphi}^{-}_{i}   
)
=
\frac{1}{\alpha_{i} }  \nabla_{I^{*}(d\alpha_{i} ) } 
\nabla_{D}^{-k} \theta_{E}  
\]
 is $s_{i} $-invariant
because 
$s_{i} (\alpha_{i} ) = -\alpha_{i}$ 
and
$s_{i} (d\alpha_{i} ) = -d\alpha_{i}$. 
Therefore
\[
{\mathbf{res} }\left( 
s_{i} (
\hat{\varphi}^{-}_{i}   
)\right)
=
s_{i}\left(
{\mathbf{res} }
(
\hat{\varphi}^{-}_{i}   
)\right)
=
{\mathbf{res} }
(
\hat{\varphi}^{-}_{i}  
) 
\]
and thus
$
s_{i} (
\hat{\varphi}^{-}_{i}   
)
=
\hat{\varphi}^{-}_{i}.   
$ 
\owari

\section{The $k$-Euler derivation } 

Let $k\in \Z_{\geq 0} $.
An \textbf{extended Catalan arrangement} $\mathit{Cat}^{k}$ 
of the type $\Phi$ is an affine arrangement defined by 
$$
\mathit{Cat}^{k}
:=
\mathit{Cat}^{k}(\Phi^{+})
:=
\{H_{\alpha, j} ~\mid~ \alpha\in \Phi^{+},\,\,\,j\in\Z, \,-k \leq j\leq k\}.
$$
Its cone 
\[
{\cal C}^{k} := {\cal C}^{k}(\Phi^{+} ) := {\bf c}\mathit{Cat}^{k} 
\]
was proved to be a free arrangement with 
\[
\exp_{0} ({\cal C}^{k} ) =(kh+d_{1}, kh+d_{2}, \dots, kh+d_{\ell}   )
\]
by Yoshinaga \cite{Y}. 
Here the integers $d_{1} , d_{2} , \dots, d_{\ell} $ are the exponents 
of the root system $\Phi$ with 
$$d_{1} \leq d_{2} \leq \dots \leq d_{\ell}. $$
Since $1=d_{1} <d_{2} $, one has 
$kh+1 < kh+d_{2} $.
Thus
we know that $D_{0} ({\cal C}^{k})_{kh+1} $ is a one-dimensional vector space. 
\begin{definition}
We say that $\eta^{k}$ is the {\bf $k$-Euler derivation} if 
  $\eta^{k}$ is a basis for the vector space $D_{0} ({\cal C}^{k})_{kh+1} $.\end{definition}
Note that the ordinary Euler derivation 
$\sum_{i=1}^{\ell} x_{i} (\partial/\partial x_{i} ) $ is a $0$-Euler derivation because 
$\mathit{Cat}^{0} = \A(\Phi)$.   The choice of a
$k$-Euler derivation has the ambiguity of nonzero constant multiples. 

\begin{prop}
\label{etakinvariant} 
A $k$-Euler derivation is $W$-invariant,
where the group $W$ acts trivially on the variable $z$. 
\end{prop}

\noindent
{\bf proof.}
Note that
${\cal C}^{k} $
is stable under the action of $W$
 unlike ${\cal S}^{k}$.
Since 
the $W$-invariant derivation
$
\nabla_{D}^{-k} \theta_{E}
$
is a basis 
for
$D (\mathcal{A}(\Phi),2k+1)_{kh+1},
$ 
we obatin
\[
D (\mathcal{A}(\Phi),2k+1)_{kh+1}
=
D (\mathcal{A}(\Phi),2k+1)_{kh+1}^{W}. 
\] 
Since the Ziegler restriction map
\[
\xymatrix{ 
D_{0} ({\mathcal C}^{k})_{kh+1} 
\ar[r]^(0.43)
{\bf {\!\!\!\!\!\!\!res}}_(0.43){\scalebox{1.5}{
\!\!\!\!\!$\widetilde{}$}}  
& 
D (\mathcal{A}(\Phi),2k+1)_{kh+1}
}
\]
is a $W$-isomorphism,
we obtain
\[
D_{0} ({\mathcal C}^{k})_{kh+1} 
=
D_{0} ({\mathcal C}^{k})_{kh+1}^{W}.  
\]
\owari
 
\medskip
We may describe a $k$-Euler derivation in terms of the SRB$+$
$
\varphi^{+}_{1},
\dots,
\varphi^{+}_{\ell}
$ 
as follows:
\begin{theorem}
\label{etak} 
The derivation
\[
\eta^{k} :=
\sum_{i=1}^{\ell}
(\alpha_{i} + kz) \varphi_{i}^{+}   
\]
is a $k$-Euler derivation. 
\end{theorem}
 
\noindent
{\bf{proof.}}
Note that ${\mathcal B}^{+}_{\Delta}$ 
is a subarrangement
of ${\mathcal C}^{k} $.
 Consider a commutative diagram
\[
\xymatrix{ 
D_{0} ({\mathcal B}^{+}_{\Delta})_{kh+1} 
\ar[r]^{\bf {\!\!\!\!\!\!\!\!\!res}}_{\scalebox{1.5}{$\!\!\!\!\! \widetilde{}$}}  
& 
D (\mathcal{A}(\Phi),{\bf z}_{\Delta}^{+})_{kh+1}
\\
D_{0} ({\mathcal C}^{k})_{kh+1} 
\ar[r]^(0.43){\bf {res}}_(0.43){\scalebox{1.5}{$\widetilde{}$}}  
\ar@{}[u]|{\scalebox{1.3}{$\bigcup$}}
& 
D(\mathcal{A}(\Phi), 2k+1)_{kh+1} 
\ar@{}[u]|{\scalebox{1.3}{$\bigcup$}}.
}
\]
Since $(\alpha_{i}+kz)\varphi^{+}_{i} \in
D_{0}({\mathcal B}^{+}_{\Delta})_{kh+1}
    $ 
for any $i$ by Theorem \ref{characterization} (1),
we have $\eta^{k} \in 
D_{0}({\mathcal B}^{+}_{\Delta})_{kh+1}
$. 
We also have
\begin{align*} 
{\mathbf{res}}\left(
\eta^{k} 
\right)
&=
{\mathbf{res}}\left(
\sum_{i=1}^{\ell}
(\alpha_{i}+kz )\varphi_{i}^{+}  
\right)\\
 &=
\sum_{i=1}^{\ell}
 \alpha_{i}  
\Xi(\partial_{\alpha^{*} _{i}})
=
\sum_{i=1}^{\ell}
 \alpha_{i}  
\nabla_{\partial_{\alpha^{*}_{i}}}
\nabla_{D}^{-k} \theta_{E}
=
\nabla_{\theta_{E}}
\nabla_{D}^{-k} \theta_{E}
\\
&=
(kh+1)   
\nabla_{D}^{-k} \theta_{E}
\in
D(\mathcal{A}(\Phi), 2k+1)_{kh+1}. 
\end{align*} 
Thus we may conclude that 
$\eta^{k} \in D_{0} ({\mathcal C}^{k})_{kh+1} 
$ 
by chasing the diagram above.
\owari

\begin{remark}
The logarithmic differential versions of Theorems \ref{characterization}
and
\ref{etak} 
can be proved similarly in the dual setup.
For the definition of the  logarithmic differential modules
$\Omega(\A(\Phi), 2k)~~(k \geq 1)$, see \cite{Z}.
For the Ziegler restriction map in the dual setup, see
\cite[Theorem 2.5]{Y05}.   
\end{remark}

The following corollary of Theorem \ref{etak} 
gives a formula for the action of a simple reflection on a 
member of the SRB+.
This formula was not covered by Theorem \ref{simplereflectionaction}.   

\begin{corollary}
\label{simplereflectionaction2}
Let $\Delta = \{\alpha_{1}, \dots, \alpha_{\ell}  \}$ 
be the set of simple roots.
Let $s_{i} $ denote the simple reflection
with respect to a simple root $\alpha_{i} $ 
with $s_{i}(\alpha_{i} )=-\alpha_{i}  
\,\,(1\leq i\leq \ell)$.  
Then 
\[
s_{i} (\varphi_{i}^{+}) 
=
\left(
\frac{\alpha_{i} + kz}{-\alpha_{i} + kz} 
\right)
{\varphi}_{i}^{+}
+
\left(
\frac{\alpha_{i}}{-\alpha_{i} + kz} 
\right)
\sum_{j\neq i}
\frac{2(\alpha_{i}, \alpha_{j}) }{(\alpha_{i}, \alpha_{i})} 
{\varphi}_{j}^{+}.
\]
\end{corollary}
 
\noindent
{\bf proof.}
Since the $k$-Euler derivation 
$
\eta^{k} =
\sum_{j=1}^{\ell}
(\alpha_{j} + kz) \varphi_{j}^{+}   
$
is $W$-invariant by Proposition \ref{etakinvariant},
we obtain
\begin{align*} 
&~~~\sum_{j=1}^{\ell}
(\alpha_{j} + kz) \varphi_{j}^{+}   
=\eta^{k} 
=s_{i} (\eta^{k}) 
=
\sum_{j=1}^{\ell}
(s_{i} (\alpha_{j}) + kz) s_{i} (\varphi_{j}^{+})\\
&=
(-\alpha_{i} + kz) s_{i} (\varphi_{i}^{+})
+
\sum_{j\neq i}
(\alpha_{j} - \frac{2(\alpha_{i}, \alpha_{j}) }{(\alpha_{i}, \alpha_{i})} 
{\alpha}_{i} 
+kz) \varphi_{j}^{+}.
\end{align*} 
Here we used the formula in
Theorem \ref{simplereflectionaction} (1).
Solving this equation for $s_{i} (\varphi_{i}^{+})$, we may verify 
the desired result. 
\owari

\medskip

\noindent
Takuro Abe\\
Department of Mechanical Engineering and Science\\
Kyoto University\\
Kyoto 606-8501, Japan.\\
email:abe.takuro.4c@kyoto-u.ac.jp
\bigskip

\noindent
Hiroaki Terao\\
Department of Mathematics\\
Hokkaido University\\
Sapporo, Hokkaido 060-0810, Japan.\\
email:terao@math.sci.hokudai.ac.jp

\end{document}